\documentclass[11pt]{article}
\usepackage{amssymb}
\usepackage{latexsym}
\begin{document}
\title{{\bf CLASSIFYING TWO-DIMENSIONAL HYPOREDUCTIVE TRIPLE ALGEBRAS}}
\author{{\bf A. Nourou Issa} \\D\'epartement de Math\'ematiques\\
Universit\'e d'Abomey-Calavi\\01 B.P. 4521 Cotonou 01, BENIN.\\E-mail: woraniss@yahoo.fr}
\date{}
\maketitle
\begin{abstract}
Two-dimensional real hyporeductive triple algebras (h.t.a.) are investigated.
A classification of such algebras is presented. As a consequence, a classification
of two-dimensional real Lie triple algebras (i.e. generalized Lie triple systems)
and two-dimensional real Bol algebras is given.
\end{abstract}
\vspace{0.5truecm}
\vspace{0.5truecm}
\par
2000 Subject Classification : Primary 17D99.
\vspace{0.5truecm}
\par
{\bf 1. Introduction.} Hyporeductive algebras were introduced by Sabinin L. V.
([7, 8]) as an infinitesimal tool for
the study of smooth hyporeductive loops which are a generalization both of smooth
Bol loops and smooth reductive loops (i.e. smooth A-loops with monoalternative
property [8]). It is shown that the fundamental vector fields of
any smooth hyporeductive loop constitute an algebra called a hyporeductive algebra
of vector fields. Further (see [1, 2, 7]) this notion has  been extended to the one of
an abstract hyporeductive triple algebra (h.t.a. for short) meaning a finite-
dimensional linear space with two binary and one ternary operations satisfying
some specific identities. It turns out that hyporeductive algebras generalize
Bol algebras and Lie triple algebras (see [10, 12] about Bol and Lie triple algebras).
\par
In this paper we consider 2-dimensional h.t.a. over the field of real numbers
(i.e. 2-dimensional real h.t.a.) and the search of clear
expressions of operations for such algebras led us to their classification (such
a classification includes the one of 2-dimensional real Lie triple algebras, Lie
triple systems and Bol algebras).
\par
Petersson H. P. [6] solved the classification problem for 2-dimensional nonassociative
algebras over arbitrary base fields, and in his approach structure constants or
multiplication tables almost never play a significant role. Underlying this classification
is the use of an isomorphism theorem and the principal Albert isotopes. One observes
that the algebras considered in that paper are binary algebras. It seems that the Petersson
approach is not applicable to the case of 2-dimensional h.t.a. because of the following
reasons. First, h.t.a. are contained in the class of tangent structures
{\it with one (or more) binary operations and a ternary operation} (this is why they
are usually called {\it binary-ternary algebras}) satisfying some compatibility
conditions. Thus, if the Petersson approach could be applied to the binary operations
of h.t.a. (under some conditions), it does not work, in general, for the ternary
operation of h.t.a. (for instance, we still do not know what is the principal Albert
isotope of a ternary operation of an algebra). Next, even for Bol algebras which are a
very particular instance of h.t.a., almost no classification results are known (the
classification of 2-dimensional real Bol algebras given in the present paper seems to
be, to our knowledge, the first one so far). Because of the nature of h.t.a., their
classification over arbitrary base fields should generalize, e.g., the one of Bol
algebras over arbitrary base fields but, unfortunately, the latter is still not
available in litterature. The other reason for considering in this paper only real
h.t.a. is related to the correspondence between h.t.a. and smooth hyporeductive loops
([9]) (this problem is solved by Kuz'min E. N. [5] for {\it real} finite-dimensional
Malcev algebras and
smooth Moufang loops and by Sabinin L. V. and Mikheev P. O. [10] for {\it real} finite-dimensional
Bol algebras and smooth
Bol loops).
\par
In Section 2 some results on hyporeductive algebras are recalled and the classification
theorem is stated. The section 3 deals with its proof (this proof gives the
classification strategy). \\
\par
{\bf 2. Background and results.} Hyporeductive algebras were originally introduced ([7, 8])
as algebras of vector fields on a smooth finite-dimensional manifold, satisfying a
specific condition. More exactly it was given the following\\
\par
{\bf Definition 1} [7, 8]. {\it A linear space $V$ of vector fields on a real $n$-dimensional
manifold $M$ with a singled out point $e$, satisfying
\begin{equation}
[X,[Y,Z]] = [X,a(Y,Z)] + r(X;Y,Z)
\end{equation}
is called a {\it hyporeductive algebra} of vector fields with determining operations
$a$ and $r$, if $dim \{X(e) : X \in V\} = n$.}
\par
Obviously $a(Y,Z)$ is a bilinear skew-symmetric operation and $r(X;Y,Z)$ a trilinear
operation on $V$, skew-symmetric in the last two variables. We called the relation (1)
the {\it hyporeductive condition} for algebras of vector fields (see [1, 2]). Considering a hyporeductive
algebra as a tangent algebra at the identity $e$ of a smooth hyporeductive loop
it is shown ([7, 9]) that a hyporeductive algebra may be viewed as an algebra with
two binary operations $a(X,Y)(e),\;T_{e}(X,Y)=[X,Y](e)$ and one ternary operation
$r(Z;X,Y)(e)$ and then, working out the Jacobi identities in the corresponding enveloping
Lie algebra, one can get the full system of identities linking the operations
$a, T_{e}, r$. A similar construction is carried out in [1, 2], where instead
of $T_{e}(X,Y)$ the operation $b(X,Y)=[X,Y](e)-a(X,Y)(e)$ is introduced (this is made
in connection with a more suitable differential geometric interpretation of a
hyporeductive algebra of vector fields and then the system of identities
mentioned above constitutes the integrability conditions of the structure equations
of the affinely connected smooth manifold associated with a local smooth hyporeductive
loop). If define the operations \\
$X \cdot Y= a(X,Y)(e)$, $X*Y=b(X,Y)$ and $\langle Z;X,Y \rangle = r(Z;X,Y)(e)$, \\
then $(T_{e}M, \cdot ,  * , \langle  ;  , \; \rangle )$ is an algebra satisfying
the system of identities mentioned above. This led us to introduce the notion of
an abstract hyporeductive triple algebra (h.t.a.): \\
\par
{\bf Definition 2} ([1, 2, 3]). Let $\cal V$ be a finite-dimensional linear space.
Assume that on ${\cal V}$ are defined two binary skew-symmetric operations \\ $"\cdot",\; "*"$
and one ternary operation $"\langle ; , \rangle "$ skew-symmetric in the last two
variables. We say that the algebra $(\cal V, \cdot , * , \langle ; , \rangle )$ is an {\it abstract}
h.t.a. if, for any $\xi, \eta, \zeta, \kappa, \chi, \theta$ in $\cal V$ the following
identities hold:
\begin{equation}
{\sigma}\;\;\; \{ \;\;{\xi}\; \cdot \;({\eta}\; \cdot \;{\zeta}) - \langle {\xi};{\eta},{\zeta} \rangle  \;\;\; \} = 0,
\end{equation}
\begin{equation}
{\sigma}\;\;\; \{ \;\;{\zeta}\;*\;({\xi}\; \cdot \;{\eta})\;\;\; \} = 0,
\end{equation}
\begin{equation}
{\sigma}\;\;\; \{ \;\; \langle {\theta};{\zeta},\;{\xi}\; \cdot \;{\eta} \rangle \;\; \} = 0,
\end{equation}
\begin{eqnarray}
{\kappa} \cdot \langle {\zeta};{\xi},{\eta} \rangle &  & - {\zeta} \cdot \langle {\kappa};{\xi},{\eta} \rangle + \langle {\zeta} \cdot {\kappa};{\xi},{\eta} \rangle  =  \nonumber \\
                                &   & \langle {\xi}*{\eta};{\zeta},{\kappa} \rangle - \langle {\zeta}*{\kappa};{\xi},{\eta} \rangle \nonumber \\
                                &   & + {\zeta}* \langle {\kappa};{\xi},{\eta} \rangle - {\kappa}* \langle {\zeta};{\xi},{\eta} \rangle  \nonumber \\
                                &   & + ({\xi}*{\eta})*({\zeta}*{\kappa}) + ({\xi}*{\eta}) \cdot ({\zeta}*{\kappa}),
\end{eqnarray}
\begin{eqnarray}
{\chi} \cdot ({\kappa} \cdot \langle {\zeta};{\xi},{\eta} \rangle &  & - {\zeta} \cdot \langle {\kappa};{\xi},{\eta} \rangle + \langle {\zeta} \cdot {\kappa};{\xi},{\eta} \rangle) \nonumber \\
                                        &  & + \langle \langle {\chi};{\xi},{\eta} \rangle ;{\zeta},{\kappa} \rangle - \langle \langle {\chi};{\zeta},{\kappa} \rangle ;{\xi},{\eta} \rangle \nonumber \\
                                        &  & + \langle {\chi};{\zeta}, \langle {\kappa};{\xi},{\eta}\rangle \rangle - \langle {\chi};{\kappa}, \langle {\zeta};{\xi},{\eta}\rangle \rangle = 0,
\end{eqnarray}
\begin{equation}
{\chi}\;*\;({\kappa} \cdot \langle {\zeta};{\xi},{\eta} \rangle - {\zeta} \cdot \langle {\kappa};{\xi},{\eta} \rangle + \langle {\zeta} \cdot {\kappa};{\xi},{\eta} \rangle ) = 0,
\end{equation}
\begin{equation}
\langle {\theta};{\chi},\;\;{\kappa} \cdot \langle {\zeta};{\xi},{\eta}\rangle - {\zeta} \cdot \langle {\kappa};{\xi},{\eta} \rangle + \langle {\zeta} \cdot {\kappa};{\xi},{\eta} \rangle \rangle = 0,
\end{equation}
\begin{eqnarray}
{\kappa} \cdot \langle {\zeta};{\xi},{\eta} \rangle &  & - {\zeta} \cdot \langle {\kappa};{\xi},{\eta} \rangle + \langle {\zeta} \cdot {\kappa};{\xi},{\eta} \rangle \nonumber \\
& & +{\eta} \cdot \langle {\xi};{\zeta},{\kappa} \rangle - {\xi} \cdot \langle {\eta};{\zeta},{\kappa} \rangle + \langle {\xi} \cdot {\eta};{\zeta},{\kappa} \rangle =0,
\end{eqnarray}
\begin{equation}
{\zeta}* \langle {\kappa};{\xi},{\eta} \rangle - {\kappa}* \langle {\zeta};{\xi},{\eta} \rangle + {\xi}* \langle {\eta};{\zeta},{\kappa} \rangle  - {\eta}* \langle {\xi};{\zeta},{\kappa} \rangle =0,
\end{equation}
\begin{eqnarray}
{\Sigma} \;\; & \{ & \;\; \langle (\langle {\xi} \cdot {\eta};{\zeta},{\kappa} \rangle +{\eta} \cdot \langle {\xi};{\zeta},{\kappa} \rangle \nonumber \\
            & - & {\xi} \cdot \langle {\eta};{\zeta},{\kappa} \rangle ); {\lambda} , {\mu} \rangle \nonumber \\
            & + & \langle {\lambda} \cdot {\mu}; \langle {\eta};{\zeta},{\kappa} \rangle , {\xi} \rangle + {\mu} \cdot \langle {\lambda}; \langle {\eta};{\zeta},{\kappa} \rangle , {\xi} \rangle \nonumber \\
            & - & {\lambda} \cdot \langle {\mu};\langle {\eta};{\zeta},{\kappa} \rangle , {\xi} \rangle \nonumber \\
            & - & (\langle {\lambda} \cdot {\mu}; \langle {\xi};{\zeta},{\kappa} \rangle , {\eta} \rangle + {\mu} \cdot \langle {\lambda}; \langle {\xi};{\zeta},{\kappa} \rangle , {\eta} \rangle \nonumber \\
            & - & {\lambda} \cdot \langle {\mu};\langle {\xi};{\zeta},{\kappa} \rangle , {\eta} \rangle ) \; \; \} =0,
\end{eqnarray}
\begin{eqnarray}
{\Sigma} & \{ & ( \langle {\mu};\langle {\eta};{\zeta},{\kappa} \rangle ,{\xi} \rangle - \langle {\mu}; \langle {\xi};{\zeta},{\kappa} \rangle ,{\eta} \rangle ) * {\lambda} \nonumber \\
         & + & ( \langle {\lambda};\langle {\xi};{\zeta},{\kappa} \rangle ,{\eta} \rangle - \langle {\lambda};\langle {\eta};{\zeta},{\kappa} \rangle ,{\xi} \rangle ) * {\mu} \;\; \} = 0,
\end{eqnarray}
\begin{eqnarray}
{\Sigma} & \{ & \langle {\theta};( \langle {\mu}; \langle {\eta};{\zeta},{\kappa} \rangle ,{\xi} \rangle - \langle {\mu};\langle {\xi};{\zeta},{\kappa} \rangle ,{\eta} \rangle ),{\lambda} \rangle \nonumber \\
         & + & \langle {\theta}; ( \langle {\lambda};\langle {\xi};{\zeta},{\kappa} \rangle ,{\eta} \rangle - \langle {\lambda};\langle {\eta};{\zeta},{\kappa} \rangle ,{\xi} \rangle ),{\mu} \rangle \} = 0,
\end{eqnarray} \\
where ${\sigma}$ denotes the sum over cyclic permutations of $\xi$, $\eta$, $\zeta$ and $\Sigma$ the one on pairs ($\xi$, $\eta$), ($\zeta$, $\kappa$), ($\lambda$, $\mu$). \\
\par
{\bf Remark 1.} The study of h.t.a. is more tractable if they are given in terms of identities as in
the definition above. For instance, we observe that if in (2)-(13) we set ${\xi} \cdot {\eta} = 0$
for any ${\xi},{\eta}$ of ${\cal V}$, then we get the defining identities of
a {\it Bol algebra} $(\cal V, * , \langle ; , \rangle )$: \\
${\sigma}\;\;\; \{ \;\;\langle {\xi};{\eta},{\zeta} \rangle \;\; \} = 0$, \\
$\langle {\xi}*{\eta};{\zeta}, {\kappa} \rangle - \langle {\zeta}* {\kappa};{\xi},{\eta} \rangle + $ \\
${\zeta}* \langle {\kappa};{\xi},{\eta} \rangle -{\kappa} * \langle {\zeta};{\xi},{\eta} \rangle + ({\xi}*{\eta})*({\zeta}* {\kappa}) = 0$, \\
$\langle \langle {\chi}; {\xi},{\eta} \rangle ;{\zeta}, {\kappa} \rangle - \langle \langle {\chi}; {\zeta},{\kappa} \rangle ;{\xi},{\eta} \rangle +$ \\
$+ \langle {\chi}; {\zeta},\langle {\kappa} ;{\xi},{\eta} \rangle \rangle - \langle {\chi}; {\kappa},\langle {\zeta};{\xi},{\eta} \rangle \rangle = 0$. \\
From the other hand, setting ${\xi}*{\eta} = 0$, we get a {\it Lie triple
algebra} (i.e. a {\it generalized Lie triple system}) $(\cal V, \cdot , \langle ; , \rangle )$: \\
${\sigma}\;\;\; \{ \;\;{\xi}\; \cdot \;({\eta}\; \cdot \;{\zeta}) - \langle {\xi};{\eta},{\zeta} \rangle \;\;\; \} = 0$, \\
${\sigma}\;\;\; \{ \;\;\langle {\theta};{\zeta},\;{\xi}\; \cdot \;{\eta} \rangle \;\; \} = 0$, \\
${\kappa} \cdot \langle {\zeta};{\xi},{\eta} \rangle - {\zeta} \cdot \langle {\kappa};{\xi},{\eta} \rangle + \langle {\zeta} \cdot {\kappa};{\xi},{\eta} \rangle  = 0$, \\
$\langle \langle {\chi};{\xi},{\eta} \rangle ;{\zeta},{\kappa} \rangle - \langle \langle {\chi};{\zeta},{\kappa} \rangle ;{\xi},{\eta} \rangle + \langle {\chi};{\zeta},\langle {\kappa};{\xi},{\eta} \rangle \rangle -$ \\
$- \langle {\chi};{\kappa}, \langle {\zeta};{\xi},{\eta} \rangle \rangle = 0$ \\
and if, moreover, we put ${\xi} \cdot {\eta} = 0$
then we obtain a {\it Lie triple system (L.t.s.)} (see Yamaguti K. [11]). Note that for
${\xi}*{\eta} = 0$ and  ${\xi} \cdot {\eta} = 0$, the identities (9)-(13) hold trivially.
\par
The question naturally arises whether there exist proper abstract h.t.a. The
answer to this problem is easier to seek among low-dimensional h.t.a. because of
the specific properties of operations $" \cdot ","*"$ and $"\langle ; , \rangle "$. Thus we are led
to the study of two-dimensional real h.t.a. that is, to find the clear expressions
of their defining operations. The following classification
theorem describes, up to isomorphisms, all 2-dimensional real h.t.a. Such a
classification includes the one of 2-dimensional real Bol algebras, Lie triple algebras
and Lie triple systems. \\
\par
{\bf Theorem}. {\it Any 2-dimensional real h.t.a. is isomorphic to one of the
h.t.a. of the following types}:\\
(I) $ u*v =0$, $u \cdot v =0$, $\langle u;u,v \rangle = eu+fv$, $ \langle v;u,v \rangle = ku-ev$, \\
(II) $ u*v =0$, $u \cdot v =au$, $\langle u;u,v \rangle = 0$, $ \langle v;u,v \rangle = ku$, \\ ($a {\neq} 0$), \\
(III) $ u*v =0$, $u \cdot v =au + bv$, $\langle u;u,v \rangle = 0$, $ \langle v;u,v \rangle = 0$, \\ ($a {\neq} 0$, $b {\neq} 0$), \\
(IV) $ u*v =0$, $u \cdot v =au+bv$, $\langle u;u,v \rangle = eu+fv$, $ \langle v;u,v \rangle = ku-ev$, \\
($a {\neq} 0$, $b {\neq} 0$, $e {\neq} 0$, $f {\neq} 0$, $k {\neq} -e$, $ af-be=0=bk+ae$), \\
(V) $ u*v =cu+dv$, $u \cdot v =0$, $\langle u;u,v \rangle = eu+fv$, $ \langle v;u,v \rangle = ku-ev$, \\
(($c,d$) ${\neq}$ ($0,0$)), \\
(VI) $ u*v =cu+dv$, $u \cdot v =au$, $\langle u;u,v \rangle = 0$, $ \langle v;u,v \rangle = ku$, \\
($a {\neq} 0$, ($c,d$) ${\neq}$ ($0,0$)), \\
(VII) $ u*v =cu+dv$, $u \cdot v =au+bv$, $\langle u;u,v \rangle = eu+fv$, $ \langle v;u,v \rangle = ku-ev$, \\
($a {\neq} 0$, $b {\neq} 0$, $e {\neq} 0$, $f {\neq} 0$, $k {\neq} 0$, ($c,d$) ${\neq}$ ($0,0$), $af-be=0=bk+ae$), \\
(VIII) $ u*v =cu+dv$, $u \cdot v =au+bv$, $\langle u;u,v \rangle = 0$, $ \langle v;u,v \rangle = 0$, \\
($a {\neq} 0$, $b {\neq} 0$, ($c,d$) ${\neq}$ ($0,0$)), \\
{\it where a, b, c, d, e, f, k are real numbers}. \\
\par
{\bf Remark 2}. The algebras of types (I), (II), (III) and (IV) are 2-dimensional
real Lie triple algebras, the zero algebra and 2-dimensional Lie triple systems
are contained in the type (I). A comprehensive classification of 2-dimensional complex
Lie triple systems is given by Yamaguti K. [11] (indeed, the type I above is just
the Lemma 5.1 of [11] when the base field is the one of real numbers); see also
Jacobson N. [4]. Note that the types (II), (IV)
are nontrivial real Lie triple algebras. The algebras of type (V) constitute nontrivial
2-dimensional real Bol algebras (see Corollary 2 below). \\
\par
{\bf Corollary 1}. {\it There exist nontrivial 2-dimensional real h.t.a. Moreover,
any such an algebra is isomorphic to an algebra of type (VI), (VII) or (VIII)}. 
\hfill ${\square}$ \\
\par
At this point we note that the example of a 2-dimensional real h.t.a. that
we gave in [3] is isomorphic to the algebra of type (VI) given by $u*v =dv$,
$u \cdot v =u$, $\langle u;u,v \rangle = 0$, $ \langle v;u,v \rangle = -u$, ($d \neq 0$).
\par
The subject of this paper is originally motivated by the need of showing concrete
nontrivial h.t.a. Besides, an affine connection space locally permitting a structure
of h.t.a. of vector fields is already described in [2] and, conversely, the structure
equations of such an affine connection space give rise to a h.t.a. structure on the
tangent space at a given point of the manifold. In relation with this, we consider
here an example of such an affine connection space with a local loop structure ([10])
with the sole condition that is given a skew-symmetric bilinear function on the space
of certain vector fields.
\par
Let $(U, \circ , e)$ be a smooth local loop so that $U$ is a sufficiently small neighborhood
of the fixed point $e$ of a real $n$-dimensional manifold $M$. We may consider on $U$
the so-called right fundamental vector fields $\{ X_{\sigma} \}$ of the loop $(U, \circ , e)$
(see, e.g., [9] and references therein), $[X_{\sigma}(x)]^{\tau}$ =
$X_{\sigma}^{\tau}(x)$, $x \in U$. Since $X_{\sigma}^{\tau}(e)$=${\delta}_{\sigma}^{\tau}$
and $e$ is a two-sided identity of $(U, \circ , e)$, it follows that $X_{1},...,X_{n}$
define a basis of vector fields linearly independent at each point of $U$ and thus $U$ is parallelizable. The Lie
bracket of two basis vector fields $X_{\alpha}$, $X_{\beta}$ is $[X_{\alpha},X_{\beta}](x)$ = $C_{\alpha \beta}^{\gamma}(x) X_{\gamma}(x)$
(observe that, in contrast of the case of left-invariant vector fields of a Lie group,
the $C_{\alpha \beta}^{\gamma}$ are functions of point [7, 9]). Now define on $U$
the (-)-connection ${\nabla}_{Z}Y = 0$ obtained from the parallelization, for any vector fields $Y$, $Z$ on $U$. Assume
that on the space of all right fundamental vector fields on $U$ is given a skew-symmetric
bilinear function $a(Y,Z)$. The torsion $T$ of the connection defined above has the
expression $T(Y,Z) =-[Y,Z]$. The vector field $[Y,Z]-a(Y,Z)$ is defined on $U$ and so
is the vector field
$[W,[Y,Z]-a(Y,Z)]$, where $W,Y,Z$ are right fundamental vector fields on $U$. Therefore,
with respect to the basis $\{ X_{1},...,X_{n} \}$, we have the representation \\
(*) $[X_{l},[X_{j},X_{k}]-a(X_{j},X_{k})]= r_{l,jk}^{i}X_{i}$, \\
which means that a structure of a h.t.a. of vector fields is locally defined (this is the
original definition of a hyporeductive algebra of vector fields [7, 8, 9]). The relation
(*) may be written as \\
(**) $({\nabla}_{l}T_{jk}^{i} - T_{ls}^{i}(T_{jk}^{s}+ a_{jk}^{s}))(x) =-r_{l,jk}^{i}(x)$ \\
for any $x \in U$, where the skew-symmetric tensor $(a_{jk}^{s})$ is defined by $a(X_{j},X_{k})= a_{jk}^{s}X_{s} $, $(T_{jk}^{s})$
is the torsion tensor of the connection $\nabla$ and ${\nabla}_{l}T_{jk}^{i}$ denotes
the covariant derivative of the function $T_{jk}^{i}$ by the vector field $X_{l}$.
Since $\{ X_{1},...,X_{n} \}$ is a parallelization, the
$r_{l,jk}^{i}$ are constants and the relation (**) means that ${\nabla}_{m}({\nabla}_{l}T_{jk}^{i} - T_{ls}^{i}(T_{jk}^{s}+ a_{jk}^{s}))=0$
at each point of $U$. The structure equations of $(U,\nabla)$ in terms of the
basis $\{ X_{1},...,X_{n} \}$ is then \\
$d{\omega}^{i}= {\frac{1}{2}} T_{jk}^{i}{\omega}^{j} \wedge {\omega}^{k} $, \\
$dT_{jk}^{i}={\nabla}_{l}T_{jk}^{i}{\omega}^{l}$. \\
The integrability conditions for these equations (at the point $e$) are precisely the
defining identities, written in terms of structure constants, of an abstract h.t.a.
and so the tangent space $T_{e}M$ is provided with a h.t.a. structure (see [2] for the general
case of an affine connection space related with a smooth local hyporeductive loop).
\par
Suppose now that dim M = 2 and choose the basis vector fields ${X_{1},X_{2}}$ such that
$[X_{1},X_{2}](e)= 2X_{1}(e) + X_{2}(e)$, $(X_{1}T_{12}^{1})(e)= 1$, $(X_{1}T_{12}^{2})
(e)= -1$, $(X_{2}T_{12}^{1})(e)= 1$, $(X_{2}T_{12}^{2})(e)= 0$. Moreover, choose the
skew-symmetric function $a(Y,Z)$ such that $a(X_{1},X_{2})(e)= X_{1}(e)+X_{2}(e)$.
Then, as indicated in the beginning of this section, we may define on $T_{e}M$ two
binary operations ${\widetilde X}_{1}\cdot {\widetilde X}_{2}= {\widetilde X}_{1}+{\widetilde X}_{2}$,
${\widetilde X}_{1}*{\widetilde X}_{2}= {\widetilde X}_{1}$ and, using (**), a ternary
operation $\langle {\widetilde X}_{1};{\widetilde X}_{1},{\widetilde X}_{2} \rangle = -{\widetilde X}_{1}+{\widetilde X}_{2}$,
$\langle {\widetilde X}_{2};{\widetilde X}_{1},{\widetilde X}_{2} \rangle = {\widetilde X}_{1}-{\widetilde X}_{2}$,
where ${\widetilde X}_{1}:= X_{1}(e) $ and ${\widetilde X}_{2} := X_{2}(e)$. It is easy
to see that the space $T_{e}M$ along with these operations constitutes a h.t.a. of type VII. \\
\par
According to the remarks above we have also the following \\
\par
{\bf Corollary 2}. {\it Any 2-dimensional real Lie triple algebra is isomorphic
to one of the algebras of the following types}: \\
(T1) $u \cdot v = 0$, $\langle u;u,v \rangle = {\alpha}u + {\beta}v$, $\langle v;u,v \rangle = {\gamma}u - {\alpha}v$, \\
(T2) $u \cdot v = u$, $\langle u;u,v \rangle = 0$, $\langle v;u,v \rangle = ku$, \\
(T3) $u \cdot v = u+v$, $\langle u;u,v \rangle = 0$, $\langle v;u,v \rangle = 0$, \\
(T4) $u \cdot v =au+bv$, $\langle u;u,v \rangle = eu+fv$, $ \langle v;u,v \rangle = ku-ev$, \\
($a {\neq} 0$, $b {\neq} 0$, $e {\neq} 0$, $f {\neq} 0$, $k {\neq} 0$, $ af-be=0=bk+ae$).
\par
{\it Any 2-dimensional real Bol algebra is isomorphic to one of the algebras
of the following types}: \\
(B1) $u*v = 0$, $\langle u;u,v \rangle = {\alpha}u + {\beta}v$, $\langle v;u,v \rangle = {\gamma}u - {\alpha}v$, \\
(B2) $u*v =cu+dv$, $\langle u;u,v \rangle = eu+fv$, $ \langle v;u,v \rangle = ku-ev$, \\
(B3) $u*v =cu+dv$, $\langle u;u,v \rangle = eu+fv$, $ \langle v;u,v \rangle = -ev$, \\
(B4) $u*v =cu+dv$, $\langle u;u,v \rangle = eu$, $ \langle v;u,v \rangle = ku-ev$, \\
(B5) $u*v =cu+dv$, $\langle u;u,v \rangle = eu$, $ \langle v;u,v \rangle = -ev$, \\
(B6) $u*v =cu+dv$, $\langle u;u,v \rangle = fv$, $ \langle v;u,v \rangle = ku$, \\
(B7) $u*v =cu+dv$, $\langle u;u,v \rangle = fv$, $ \langle v;u,v \rangle = 0$, \\
(B8) $u*v =cu+dv$, $\langle u;u,v \rangle = 0$, $ \langle v;u,v \rangle = 0$, \\
(B9) $u*v =cu$, $\langle u;u,v \rangle = eu+fv$, $ \langle v;u,v \rangle = ku-ev$, \\
(B10) $u*v =cu$, $\langle u;u,v \rangle = eu+fv$, $ \langle v;u,v \rangle = -ev$, \\
(B11) $u*v =cu$, $\langle u;u,v \rangle = eu$, $ \langle v;u,v \rangle = ku-ev$, \\
(B12) $u*v =cu$, $\langle u;u,v \rangle = eu$, $ \langle v;u,v \rangle =-ev$, \\
(B13) $u*v =cu$, $\langle u;u,v \rangle = fv$, $ \langle v;u,v \rangle = ku$, \\
(B14) $u*v =cu$, $\langle u;u,v \rangle = fv$, $ \langle v;u,v \rangle = 0$, \\
(B15) $u*v =cu$, $\langle u;u,v \rangle = 0$, $ \langle v;u,v \rangle = ku$, \\
(B16) $u*v =cu$, $\langle u;u,v \rangle = 0$, $ \langle v;u,v \rangle = 0$, \\
(B17) $u*v =dv$, $\langle u;u,v \rangle = eu+fv$, $ \langle v;u,v \rangle =-ev$, \\
(B18) $u*v =dv$, $\langle u;u,v \rangle = eu$, $ \langle v;u,v \rangle = ku-ev$, \\
(B19) $u*v =dv$, $\langle u;u,v \rangle = eu$, $ \langle v;u,v \rangle = -ev$, \\
{\it where $c {\neq} 0$, $d {\neq} 0$, $e {\neq} 0$, $f {\neq} 0$, $k {\neq} 0$}. \hfill ${\square}$ \\
\par
The types (B2)-(B19) constitute just the developed form of the type (V). \\
\par
{\bf 3. Proof of the theorem.} First we shall prove the following \\
{\bf Lemma}. {\it If $\{x_{1},x_{2}\}$ is a basis of a 2-dimensional real h.t.a. ${\cal V}$,
then the identities (2)-(13) of abstract h.t.a. have the following form:
\begin{eqnarray}
J(x_{1},x_{2})&  &- \langle x_{1} \cdot x_{2};x_{1},x_{2} \rangle  \nonumber \\
              &  & +x_{1} \cdot \langle x_{2};x_{1},x_{2} \rangle-x_{2} \cdot \langle x_{1};x_{1},x_{2} \rangle =0,
\end{eqnarray}
\begin{eqnarray}
x_{i} \cdot J(x_{1},x_{2})&  &- \langle x_{i};x_{1},\langle x_{2};x_{1},x_{2} \rangle \rangle  \nonumber \\
                    &  &+ \langle x_{i};x_{2}, \langle x_{1};x_{1},x_{2} \rangle \rangle =0,
\end{eqnarray}
\begin{equation}
x_{i}*J(x_{1},x_{2})=0,
\end{equation}
\begin{equation}
\langle x_{j};x_{i},J(x_{1},x_{2}) \rangle =0,
\end{equation}
\begin{equation}
\langle x_{1} \cdot x_{2};x_{1},x_{2} \rangle - x_{1} \cdot \langle x_{2};x_{1},x_{2} \rangle + x_{2} \cdot \langle x_{1};x_{1},x_{2} \rangle =0,
\end{equation}
\begin{equation}
J(x_{1},x_{2})=0,
\end{equation}
where $J(x_{1},x_{2})= x_{1}*\langle x_{2};x_{1},x_{2} \rangle -x_{2}*\langle x_{1};x_{1},x_{2} \rangle $ and $i,j=1,2$.}\\
{\it Proof}. With respect to the basis $\{x_{1},x_{2}\}$, (2),(3) and (4) are
clearly satisfied trivially. Next the left-hand side of (5) now reads $x_{i} \cdot \langle x_{j};x_{1},x_{2} \rangle$ \\ $-x_{j} \cdot \langle x_{i};x_{1},x_{2} \rangle +\langle x_{j} \cdot x_{i};x_{1},x_{2} \rangle $
while the right-hand side reads \\ $\langle x_{1}*x_{2};x_{j},x_{i} \rangle - \langle x_{j}*x_{i};x_{1},x_{2} \rangle +x_{j}*\langle x_{i};x_{1},x_{2} \rangle - \\ -x_{i}*\langle x_{j};x_{1},x_{2} \rangle +(x_{1}*x_{2})*(x_{j}*x_{i})+(x_{1}*x_{2}) \cdot (x_{j}*x_{i})$,
with $i,j=1,2$. Furthermore, because of the skew-symmetry of operations
$" \cdot ", "*", \\ "\langle ; , \rangle "$ one observes that the identity (5) gets the form \\
$x_{1}*\langle x_{2};x_{1},x_{2} \rangle -x_{2}*\langle x_{1};x_{1},x_{2} \rangle =\langle x_{1} \cdot x_{2};x_{1},x_{2} \rangle -x_{1} \cdot \langle x_{2};x_{1},x_{2} \rangle +x_{2} \cdot \langle x_{1};x_{1},x_{2} \rangle $,
so we obtain (14).\\
In view of (14), the identities (7) and (8) are straightforwardly transformed into
(16) and (17) respectively.\\
Finally, and again with (14) in mind, we work the identity (6) as follows: we
replace ${\xi},{\eta},{\zeta},{\kappa}, {\chi}$ by $x_{1},x_{2},x_{k},x_{j},x_{i}$
respectively where $i,j,k=1,2$ and then by (14), we see that (6) gets
the form \\ $x_{i} \cdot (x_{1}* \langle x_{2};x_{1},x_{2} \rangle -x_{2}* \langle x_{1};x_{1},x_{2} \rangle ) + \langle x_{i};x_{2}, \langle x_{1};x_{1},x_{2} \rangle \rangle \\ - \langle x_{i};x_{1}, \langle x_{2};x_{1},x_{2}\rangle \rangle =0$
that is, we get (15). The equation (10) implies (19) and hence (18),
in view of (14). The equalities (11)-(13) hold trivially in view of (19) and
(15). \hfill ${\square}$ \\
\par
One observes that (18) and (19) are actually equivalent and accordingly the
system (14)-(19) takes a simpler form (see the theorem's proof below). We keep
the system (14)-(19) as above in order to follow the step by step transformation
of the system (2)-(13) in the 2-dimensional case. We now turn to the proof of
the theorem. \\
\par
{\it Proof of the theorem}.
\par
Let $(\cal V, \cdot , * , \langle ; , \rangle )$ be a 2-dimensional real h.t.a. with basis $ \{ u,v \} $.
Put  $u \cdot v =au+bv$, $u*v =cu+dv$, $\langle u;u,v \rangle = eu+fv$, $ \langle v;u,v \rangle = ku+lv$. Then,
by the lemma, the identities (2)-(13) reduce to (14)-(19) and a careful
reading of these identities reveals that the expression $ N = \langle u \cdot v;u,v \rangle $ \\ $+ v \cdot \langle u;u,v \rangle -u \cdot \langle v;u,v \rangle $
can be conclusive for the study of h.t.a. (at least in the 2-dimensional case).
Therefore we shall discuss the case $N=0$ (see(18)).
\par
Thus $ 0=N = \langle u \cdot v;u,v \rangle + v \cdot \langle u;u,v \rangle -u \cdot \langle v;u,v \rangle = $ \\
$\langle au+bv;u,v \rangle + v \cdot (eu+fv)-u \cdot (ku+lv) = (bk-al)u+(af-be)v$ implies
\begin{eqnarray}
bk-al & = & 0, \nonumber \\
af-be & = & 0.
\end{eqnarray}
Discussing the solutions of the system (20), we see that the following essential
situations occur (any other situation is either one of those enumerated below
or is included in some of them): \\
(3.1) $a {\neq} 0$, $b {\neq} 0$, $e {\neq} 0$, $f {\neq} 0$, $k {\neq} 0$, $l {\neq} 0$, \\
(3.2) $a {\neq} 0$, $b {\neq} 0$, $e {\neq} 0$, $f {\neq} 0$, $k = 0$, $l = 0$, \\
(3.3) $a {\neq} 0$, $b {\neq} 0$, $e = 0$, $f = 0$, $k {\neq} 0$, $l {\neq} 0$, \\
(3.4) $a {\neq} 0$, $b {\neq} 0$, $e = 0$, $f = 0$, $k = 0$, $l = 0$, \\
(3.5) $a {\neq} 0$, $b = 0$, $e \; any $, $f = 0$, $k = 0$, $l = 0$, \\
(3.6) $a = 0$, $b {\neq} 0$, $e = 0$, $f \; any$, $k {\neq} 0$, $l \; any$, \\
(3.7) $a = 0$, $b = 0$, $e \; any$, $f \; any$, $k \; any$, $l \; any$.
\par
Now each of the cases (3.1)-(3.7) must be discussed in connection with the
identities (14)-(19). We observe that, with the condition $N=0$ (i.e. (18))
only (14) and (15) are of interest here.
\par
The identity (14) implies $ u* \langle v;u,v \rangle -v* \langle u;u,v \rangle = 0$ which means that
\begin{equation}
(e+l)u*v =0.
\end{equation}
The equation (21) yields the following cases: \\
(3.8) $u*v=0$ and $l \neq -e$, \\
(3.9) $u*v=0$ and $l = -e$, \\
(3.10) $u*v=cu+dv$ (($c,d$) $\neq$ (0,0)) and $l = -e$. \\
Then considering each of the cases (3.1)-(3.7) in connection with the conditions
(3.8)-(3.10), we are led to the following types of algebras: \\
(A1) $u*v=0$, $u \cdot v =au+bv$, $\langle u;u,v \rangle = eu+fv$, $ \langle v;u,v \rangle = ku+lv$, \\
($a {\neq} 0$, $b {\neq} 0$, $e {\neq} 0$, $f {\neq} 0$, $k {\neq} 0$,
$l {\neq} 0$, $l {\neq} -e$, $ af-be=0=bk-al$), \\
(A2) $u*v=0$, $u \cdot v =au+bv$, $\langle u;u,v \rangle = eu+fv$, $ \langle v;u,v \rangle = ku-ev$, \\
($a {\neq} 0$, $b {\neq} 0$, $e {\neq} 0$, $f {\neq} 0$, $k {\neq} 0$,
$ af-be=0=bk+ae$), \\
(A3) $u*v=cu+dv$, $u \cdot v =au+bv$, $\langle u;u,v \rangle = eu+fv$, $ \langle v;u,v  \rangle = ku-ev$, \\
($a {\neq} 0$, $b {\neq} 0$,  ($c,d$) $\neq$ (0,0), $e {\neq} 0$, $f {\neq} 0$, $k {\neq} 0$,
$ af-be=0=bk+ae$), \\
(A4) $u*v=0$, $u \cdot v =au+bv$, $\langle u;u,v \rangle = eu+fv$, $ \langle v;u,v \rangle = 0$, \\
($a {\neq} 0$, $b {\neq} 0$, $e {\neq} 0$, $f {\neq} 0$, $k {\neq} 0$,
$ af-be=0$), \\
(A5) $u*v=0$, $u \cdot v =au+bv$, $\langle u;u,v \rangle = 0$, $ \langle v;u,v \rangle = ku+lv$, \\
($a {\neq} 0$, $b {\neq} 0$, $k {\neq} 0$, $l \neq 0$, $ bk-al =0$), \\
(A6) $u*v=0$, $u \cdot v =au+bv$, $\langle u;u,v \rangle = 0$, $ \langle v;u,v \rangle = 0$, \\
($a {\neq} 0$, $b {\neq} 0$), \\
(A7) $u*v=cu+dv$, $u \cdot v =au+bv$, $\langle u;u,v \rangle = 0$, $ \langle v;u,v \rangle =0$, \\
($a {\neq} 0$, $b {\neq} 0$,  ($c,d$) $\neq$ (0,0)), \\
(A8) $u*v=cu+dv$, $u \cdot v =au$, $\langle u;u,v \rangle = 0$, $ \langle v;u,v \rangle = ku$, \\
($a {\neq} 0$, ($c,d$) $\neq$ (0,0)), \\
(A9) $u*v=0$, $u \cdot v =au$, $\langle u;u,v \rangle = 0$, $ \langle v;u,v \rangle = ku$, \\
($a {\neq} 0$), \\
(A10) $u*v=0$, $u \cdot v =au$, $\langle u;u,v \rangle = eu$, $ \langle v;u,v \rangle = ku$, \\
($a {\neq} 0$, $e {\neq} 0$), \\
(A11) $u*v=cu+dv$, $u \cdot v =bv$, $\langle u;u,v \rangle = fv$, $ \langle v;u,v \rangle = 0$, \\
($b {\neq} 0$, ($c,d$) $\neq $ (0,0)), \\
(A12) $u*v=0$, $u \cdot v =bv$, $\langle u;u,v \rangle = fv$, $ \langle v;u,v \rangle = 0$, \\
($b {\neq} 0$), \\
(A13) $u*v=0$, $u \cdot v =bv$, $\langle u;u,v \rangle = fv$, $ \langle v;u,v \rangle = lv$, \\
($b {\neq} 0$, $l {\neq} 0$), \\
(A14) $u*v=cu+dv$, $u \cdot v =0$, $\langle u;u,v \rangle = eu+fv$, $ \langle v;u,v \rangle = ku-ev$, \\
(($c,d$) $\neq$ (0,0)), \\
(A15) $u*v=0$, $u \cdot v =0$, $\langle u;u,v \rangle = eu+fv$, $ \langle v;u,v \rangle = ku-ev$, \\
(A16) $u*v=0$, $u \cdot v =0$, $\langle u;u,v \rangle = eu+fv$, $ \langle v;u,v \rangle = ku+lv$, \\
($l \neq -e$).
\par
Furthermore, the identity (15) implies $\langle x;\langle u;u,v \rangle ,v \rangle + \\ \langle x;u,\langle v;u,v \rangle \rangle = 0$
with $x=u$ or $v$, that is,
\begin{equation}
(e+l)\langle x;u,v \rangle = 0,
\end{equation}
$x=u$ or $v$. The equation (22) gives the following cases: \\
(3.11) $l=-e$ and $\langle x;u,v \rangle \neq 0$, $x=u$ or $v$, \\
(3.12) $l=-e$ and $\langle u;u,v \rangle \neq 0$, $\langle v;u,v \rangle =0$, \\
(3.13) $l=-e$ and $\langle u;u,v \rangle = 0$, $\langle v;u,v \rangle \neq 0$, \\
(3.14) $l=-e$ and $\langle x;u,v \rangle = 0$, $x=u$ or $v$, \\
(3.15) $l \neq -e$ and $\langle x;u,v \rangle = 0$, $x=u$ or $v$. \\
Therefore, in view of constraints (3.11)-(3.15), the algebras of types (A1), (A4),
(A5), (A10), (A13), and (A16) must be cancelled out. Now, observing that an algebra
of type (A8) is isomorphic to the one of type (A11) and an algebra of type
(A9) is isomorphic to the one of type (A12), we are left with the algebras
of types (A2), (A3), (A6), (A7), (A8), (A9), (A14), (A15). These are precisely
the ones enumerated in our classification theorem.  \hfill ${\square}$ \\
\par
{\sc {\bf Acknowledgement.}} The author wishes to thank the referees for their comments
and suggestions that help to improve the initial version of the present paper. Thank
also goes to Petersson H. P. for making available his relevant article.

\end{document}